\documentclass[12pt]{amsart}

\usepackage[centertags]{amsmath}
\usepackage{amsthm,amsfonts}

\usepackage{hyperref}
\hypersetup{
	colorlinks=true,
	urlcolor=blue,
	citecolor=blue}

\usepackage{amssymb}
\usepackage{epsfig}
\usepackage{amssymb,latexsym}
\usepackage[dvipsnames]{xcolor}

\begin{document}

\theoremstyle{plain}
\newtheorem{theorem}{Theorem}[section]
\newtheorem{lemma}[theorem]{Lemma}
\newtheorem{corollary}[theorem]{Corollary}
\newtheorem{proposition}[theorem]{Proposition}
\newtheorem{question}[theorem]{Question}
\theoremstyle{definition}
\newtheorem{prop}{Proposition}[theorem]
\newtheorem{obs}{Observation}[theorem]
\newtheorem{rem}{Remark}[theorem]
\newtheorem{cor}{Corollary}[theorem]
\newtheorem{ex}{Example}[theorem] 
\newtheorem{definition}{Definition}[theorem]

\newcommand{\binomial}[2]{\left(\begin{array}{c}#1\\#2\end{array}\right)}
\newcommand{\zar}{{\rm zar}}
\newcommand{\an}{{\rm an}}
\newcommand{\red}{{\rm red}}
\newcommand{\codim}{{\rm codim}}
\newcommand{\rank}{{\rm rank}}
\newcommand{\Pic}{{\rm Pic}}
\newcommand{\Div}{{\rm Div}}
\newcommand{\Hom}{{\rm Hom}}
\newcommand{\im}{{\rm im}}
\newcommand{\Spec}{{\rm Spec}}
\newcommand{\sing}{{\rm sing}}
\newcommand{\reg}{{\rm reg}}
\newcommand{\Char}{{\rm char}}
\newcommand{\Tr}{{\rm Tr}}
\newcommand{\res}{{\rm res}}
\newcommand{\tr}{{\rm tr}}
\newcommand{\supp}{{\rm supp}}
\newcommand{\Gal}{{\rm Gal}}
\newcommand{\Min}{{\rm Min \ }}
\newcommand{\Max}{{\rm Max \ }}
\newcommand{\Span}{{\rm Span  }}

\newcommand{\Frob}{{\rm Frob}}
\newcommand{\lcm}{{\rm lcm}}


\long\def\symbolfootnote[#1]#2{\begingroup%
\def\thefootnote{\fnsymbol{footnote}}\footnote[#1]{#2}\endgroup}

\newcommand{\soplus}[1]{\stackrel{#1}{\oplus}}
\newcommand{\dlog}{{\rm dlog}\,}    
\newcommand{\limdir}[1]{{\displaystyle{\mathop{\rm
lim}_{\buildrel\longrightarrow\over{#1}}}}\,}
\newcommand{\liminv}[1]{{\displaystyle{\mathop{\rm
lim}_{\buildrel\longleftarrow\over{#1}}}}\,}
\newcommand{\boxtensor}{{\Box\kern-9.03pt\raise1.42pt\hbox{$\times$}}}
\newcommand{\sext}{\mbox{${\mathcal E}xt\,$}}
\newcommand{\shom}{\mbox{${\mathcal H}om\,$}}
\newcommand{\coker}{{\rm coker}\,}
\renewcommand{\iff}{\mbox{ $\Longleftrightarrow$ }}
\newcommand{\onto}{\mbox{$\,\>>>\hspace{-.5cm}\to\hspace{.15cm}$}}

\newenvironment{pf}{\noindent\textbf{Proof.}\quad}{\hfill{$\Box$}}

\newcommand{\sA}{{\mathcal A}}
\newcommand{\sB}{{\mathcal B}}
\newcommand{\sC}{{\mathcal C}}
\newcommand{\sD}{{\mathcal D}}
\newcommand{\sE}{{\mathcal E}}
\newcommand{\sF}{{\mathcal F}}
\newcommand{\sG}{{\mathcal G}}
\newcommand{\sH}{{\mathcal H}}
\newcommand{\sI}{{\mathcal I}}
\newcommand{\sJ}{{\mathcal J}}
\newcommand{\sK}{{\mathcal K}}
\newcommand{\sL}{{\mathcal L}}
\newcommand{\sM}{{\mathcal M}}
\newcommand{\sN}{{\mathcal N}}
\newcommand{\sO}{{\mathcal O}}
\newcommand{\sP}{{\mathcal P}}
\newcommand{\sQ}{{\mathcal Q}}
\newcommand{\sR}{{\mathcal R}}
\newcommand{\sS}{{\mathcal S}}
\newcommand{\sT}{{\mathcal T}}
\newcommand{\sU}{{\mathcal U}}
\newcommand{\sV}{{\mathcal V}}
\newcommand{\sW}{{\mathcal W}}
\newcommand{\sX}{{\mathcal X}}
\newcommand{\sY}{{\mathcal Y}}
\newcommand{\sZ}{{\mathcal Z}}

\newcommand{\A}{{\mathbb A}}
\newcommand{\B}{{\mathbb B}}
\newcommand{\C}{{\mathbb C}}
\newcommand{\D}{{\mathbb D}}
\newcommand{\E}{{\mathbb E}}
\newcommand{\F}{{\mathbb F}}
\newcommand{\G}{{\mathbb G}}
\newcommand{\HH}{{\mathbb H}}
\newcommand{\I}{{\mathbb I}}
\newcommand{\J}{{\mathbb J}}
\newcommand{\M}{{\mathbb M}}
\newcommand{\N}{{\mathbb N}}
\renewcommand{\P}{{\mathbb P}}
\newcommand{\Q}{{\mathbb Q}}
\newcommand{\T}{{\mathbb T}}
\newcommand{\U}{{\mathbb U}}
\newcommand{\V}{{\mathbb V}}
\newcommand{\W}{{\mathbb W}}
\newcommand{\X}{{\mathbb X}}
\newcommand{\Y}{{\mathbb Y}}
\newcommand{\Z}{{\mathbb Z}}


\newcommand{\Fqm}{\mathbb{F}_{q^m}}
\newcommand{\Fq}{\mathbb{F}_q}
\newcommand{\Fp}{\mathbb{F}_p}
\newcommand{\Fpl}{\mathbb{F}_{p^l}}
\newcommand{\fqn}{\mathbb{F}_q^n}
\newcommand{\be}{\begin{eqnarray}}
\newcommand{\ee}{\end{eqnarray}}
\newcommand{\nn}{{\nonumber}}
\newcommand{\dd}{\displaystyle}
\newcommand{\ra}{\rightarrow}
\newcommand{\bigmid}[1][12]{\mathrel{\left| \rule{0pt}{#1pt}\right.}}
\newcommand{\cl}{${\rm \ell}$}
\newcommand{\clp}{${\rm \ell^\prime}$}
\title{2-Colored Rogers-Ramanujan Partition Identities}
\author{MOHAMMAD ZADEHDABBAGH,  \\
}

\maketitle
\begin{center}
\fontsize{9pt}{9pt}\selectfont
Department of Mathematics, Faculty of Engineering and Natural Sciences, Sabanci University, Istanbul, Turkey\\
E-mail: mzadehdabbagh@sababnciuniv.edu\\
\end{center}
\begin{abstract}
In this paper, we combined two types of partitions and introduced 2-colored Rogers-Ramanujan partitions. By finding some functional equations and using a constructive method, some identities have been found. Some Overpartition identities coincide with our findings. A correspondence between colored partitions and overpartitions is provided.
\end{abstract}

\section{Introduction}
\label{Sec:1}
Agarwal and Andrews \cite{RRIN} gave some identities for Rogers-Ramanujan type partitions with n copies of n, which correspond to colored partitions. Sandon and Zanello \cite{WCPI} proved identities on the colored case of some partition types. This gave us the idea of combining colored partitions and Roger-Ramanujan type partitions and apply a constructive method to find identities on those partitions.

Rogers-Ramanujan identities can be interpreted using partitions in which the difference between every two consecutive parts is at least 2. We combined these types of partitions and colored ones to have the following partition type,

\begin{definition}
	A 2-colored Rogers-Ramanujan partition of $n$ consists of two separated parts, each of the same color, and the difference between every two consecutive parts of the same color is at least two, moreover, parts in different colors do not intersect.
\end{definition}

As an example, all 2-colored Rogers-Ramanujan partitions of 6 are

$$6, \color{red}6 \color{black}, 5 + 1, \color{red}5\color{black}+1, 5 + \color{red}1\color{black}, \color{red}5\color{black}+\color{red}1\color{black}, 4 + 2, \color{red}4\color{black}+2, 4 + \color{red}2\color{black}, \color{red}4\color{black}+\color{red}2\color{black}, 3 + \color{red}2\color{black}+1, \color{red}3\color{black}+2 + \color{red}1$$

A constructive way from Kur\c{s}ung\"{o}z's papers \cite{AndSty} and \cite{BreSty} applied on the generating functions for these types of partitions to find identities on 2-colored Rogers-Ramanujan partitions.

For this purpose, functional equations relating to generating functions for given partition types have been found. By these functional equations, in two main steps, we can find the generating function as an infinite sum. Finally, by some transformation such as Jacobi's triple products \cite{tToP}, we can find the divisibility part of identities.

\begin{theorem}[Jacobi's triple product]\label{JTP}
	For $z \neq 0$, $|q|<1$,
	$$\sum_{n=-\infty}^{\infty}z^nq^{n^2}=\prod_{n=0}^{\infty}(1-q^{2n+2})(1+zq^{2n+1})(1+z^{-1}q^{2n+1}).$$
\end{theorem}

By the method and the transformation mentioned above we reached the following identity

\begin{theorem}
	Let $R(n)$ denote the number of 2-colored Rogers-Ramanujan partitions of $n$, then
	$$R(n)=\frac{(-q)_\infty(q^2,q^2,q^4;q^4)_\infty}{(q)_\infty}.$$
\end{theorem}

We observed that our results are identical to identities for special cases of overpartitions in \cite{RRGO}, that identity is

\begin{theorem}
	For $k \geq a \geq 1$, let $D_{k,a}(n)$ denote the number of overpartitions of $n$ of the form $d_1 + d_2 + \cdots + d_s$, such that $1$ can occur as a non-overlined part at most $a-1$ times, and $d_j - d_{j+k-1} \geq 1$ if $d_j$ is overlined and $d_j - d_{j+k-1} \geq 2$ otherwise. For $k > i \geq 1$, let $C_{k,i}(n)$ denote the number of overpartitions of $n$ whose non-overlined parts are not congruent to $0,\pm i$ modulo $2k$ and let $C_{k,k}(n)$ denote the number of overpartitions of $n$ with parts not divisible by $k$. Then $C_{k,i}(n) = D_{k,i}(n)$.
\end{theorem}

In the first section, we will go over the 2-colored Rogers-Ramanujan partition type, accordingly, we will find two functional equations, and then constructively, mentioned above, we will find a partition identity. In the third section, a correspondence between our identities and one for overpartitions is given. Finally, possible future works are introduced in the last section.

\section{Colored Rogers-Ramanujan Partitions}
\label{Sec:2}
Throughout this paper, for $|q| < 1$,
$$(a;q)_n=\prod_{i=0}{n-1}(1-aq^i),$$
$$(a;q)_\infty=\lim_{n\rightarrow \infty}(a;q)_n,$$
$$(a_1, a_2,\cdots,a_k; q)_n := (a_1; q)_n(a_2; q)_n \cdots (a_k; q)_n.$$

According to 2-colored Rogers-Ramanujan partitions, the following definition is given.
\begin{definition}\label{R12}
	For $1 \leq j \leq 2$, let $R_j(x)$ be the generating function of 2-colored Rogers-Ramanujan partitions with smallest part greater than or equal to $j$.
\end{definition}

With respect to these definitions, one can find the following functional equations relating $R_1(x)$ and $R_2(x)$,

\begin{equation}\label{con.eq1}
	R_1(x) - R_2(x) = xq R_1(xq) + xq R_2(xq)
\end{equation}

and

\begin{equation}\label{ck.eq1}
	R_2(x) = R_1(xq)
\end{equation}

Equation \eqref{ck.eq1} is clear, as shifting every part of $R_1$ by 1 unit it will change to $R_2$. The proof of equation \eqref{con.eq1} is as follows.

Let

$$R_i(x)=\sum_{m\geq 0}\sum_{n\geq 0}r_i(m,n)x^mq^n \;\; ;\;\; i = 1, 2$$

be the generating function for the types that have been mentioned above, where m is referring to the number of parts in partitions.

Let $\pi$ be a 2-Colored Rogers-Ramanujan partition of $n$ with $m$ parts, all possible ways such that the smallest part $\geq 1$ will be counted by $r_1(m,n)$, if all partitions that the smallest part is $\geq 2$ ($r_2(m, n)$) have been removed, then all remaining partitions are of 2-Colored Rogers-Ramanujan type with the smallest part exactly be 1 ($r_1(m, n) - r_2(m, n)$).

Now, let's count all partitions with the smallest part 1 in another way, if 1 is removed from all partitions, then there will be two cases

\begin{itemize}
	\item[(i)] The smallest part is $\geq 2$ with different color than 1, so one can subtract 1 from each part, the enumeration of these partitions is by $r_1(m - 1, n - m)$.
	\item[(ii)] The smallest part is $\geq 3$ with the same color as 1, if 1 is subtracted from each part, the enumeration of these partitions is $r_2(m-1, m-n)$, note that a part 2 is not possible here.
\end{itemize}

So

$$r_1(m, n) - r_2(m, n) = r_1(m - 1, n - m) + r_2(m - 1, n - m)$$

Multiplying all terms by $x^m q^n$ and taking the summation over $m$ and $n$ for all terms, $m,n \geq 0$ and both integers, we will have

\begin{equation*}
	\begin{split}
		\sum_{m\geq 0}\sum_{n\geq 0}r_1(m, n)x^mq^n - \sum_{m\geq 0}\sum_{n\geq 0}r_2(m, n)x^mq^n =\\
		\sum_{m\geq 0}\sum_{n\geq 0}r_1(m - 1, n - m)x^mq^n+\sum_{m\geq 0}\sum_{n\geq 0}r_2(m - 1, n - m)x^mq^n
	\end{split}
\end{equation*}

By changing $m - 1$ to $m$ and $n - m$ to $n - m + 1$ on the right hand side of this equation, we will have

\begin{equation*}
	\begin{split}
		\sum_{m\geq 0}\sum_{n\geq 0}r1_(m, n)x^mq^n -\sum_{m\geq 0}\sum_{n\geq 0}r_2(m, n)x^mq^n = \\
		\sum_{m\geq 0}\sum_{n\geq 0}r_1(m, n)x^{m+1}q^{n+m+1}+\sum_{m\geq 0}\sum_{n\geq 0}r_2(m, n)x^{m+1}q^{n+m+1}
	\end{split}
\end{equation*}

so, one will get the functional equation as follows 

$$R_1(x) - R_2(x) = xqR_1(xq) + xqR_2(xq)$$

Using steps described on \cite{AndSty} and \cite{BreSty}, with straightforward  but long computations which we skipped here, by

$$R_i(x)=\sum_{n\geq 0} \alpha_n(x)q^{nA_i}+\beta_n(x)x^{B_i}q^{C_i}q^{nD_i}, \;\; i=1,2$$
we can find $\alpha_n$ and $\beta_n$ in terms of $\alpha_0$,
$$\alpha_n(x) = \tilde{\alpha}_0(xq^{2n}) \frac{x^{2n} q^{n(2n+1)} q^{-nA_2} (-1; q^E)_{n}(-x^F q^G q^{F-H} q^{nH}; q^{2F -H})_\infty}{(q^E; q^E)_n(x^F q^G q^{F-H} q^{nH}; q^{2F-H})_\infty}$$
and
$$\beta_n(x) = -\tilde{\alpha}_0(xq^{2n+1})\frac{x^{2n+1}q^{(n+1)(2n+1)}x^{-B_2} q^{-C_2} q^{-nD_2} (-1; q^E)_{n+1}(-x^Fq^Gq^{2F-H}q^{nH}; q^{2F -H})_\infty}{(q^E; q^E)_n(x^F q^G q^{nH}; q^{2F-H})_\infty}$$
where
$$\tilde{\alpha}_0(xq^{2n}) = \alpha_0(xq^{2n})\frac{((xq^{2n})^F q^G q^{F -H}; q^{2F-H})_\infty}{(-(xq^{2n})^F q^G q^{F-H}; q^{2F -H})_\infty}.$$

Then, by equation \eqref{ck.eq1} and considering some assumptions for equations to be consistent, we can find $E$, $F$, $G$ and $H$, in this case $F = G = H = 1$ and $E = -1$.

Putting them in the generating functions for $R_i(x)$, being constant of $\tilde{\alpha}_0$ with respect to $x$, and applying $x=1$, we will have
$$R_1(1)=\sum_{n\geq 0}\frac{(-1)^nq^{n(2n+1)}(-1;q)_n(-q^{n+1};q)_\infty}{(q)_n(q^{n+1};q)_\infty}-\frac{(-1)^nq^{(n+1)(2n+2)}(-1;q)_{n+1}(-q^{n+2};q)_\infty}{(q)_n(q^{n+1};q)_\infty}$$
and
$$R_2(1)=\sum_{n\geq 0}\frac{(-1)^nq^{n(2n+2)}(-1;q)_n(-q^{n+1};q)_\infty}{(q)_n(q^{n+1};q)_\infty}-\frac{(-1)^nq^{(n+1)(2n+1)}(-1;q)_{n+1}(-q^{n+2};q)_\infty}{(q)_n(q^{n+1};q)_\infty}.$$
they can be rewritten as follows,
$$R_1(1)=2\frac{(-q)_\infty}{(q)_\infty}\sum_{n\geq 0}(-1)^nq^{n(2n+1)}(\frac{1}{1+q^n}-\frac{q^{3n+2}}{1+q^{n+1}})$$
and
$$R_2(1)=2\frac{(-q)_\infty}{(q)_\infty}\sum_{n\geq 0}(-1)^nq^{n(2n+2)}(\frac{1}{1+q^n}-\frac{q^{n+1}}{1+q^{n+1}}).$$

So,

\begin{equation*}
	\begin{split}
		R_1(1)=2\frac{(-q)_\infty}{(q)_\infty}(\sum_{n\geq 0}\frac{(-1)^nq^{n(2n+1)}}{1+q^n}-\sum_{n\geq 0}\frac{(-1)^nq^{n(2n+1)}q^{3n+2}}{1+q^{n+1}})\\
		=\frac{(-q)_\infty}{(q)_\infty}(1+2(\sum_{n\geq 1}\frac{(-1)^nq^{n(2n+1)}}{1+q^n}-\sum_{n\geq 1}\frac{(-1)^{n-1}q^{2n^2}}{1+q^{n}}))\\
		=\frac{(-q)_\infty}{(q)_\infty}(1+2\sum_{n\geq 1}(-1)^nq^{2n^2})=\frac{(-q)_\infty}{(q)_\infty}\sum_{n=-\infty}^{\infty}(-1)^nq^{2n^2}
	\end{split}
\end{equation*}

By theorem \ref{JTP} for $z=-1$ and $q^2$, for 2-colored Rogers-Ramanujan defined in \ref{R12} the following identity holds

$$R_1(1)=\frac{(-q)_\infty(q^2,q^2,q^4;q^4)_\infty}{(q)_\infty}$$

Moreover, the coefficients in the Taylor series of $R_2(1)$ coincides with the number of partitions for 2-colored Rogers-Ramanujan partitions with parts more than $1$,

$$1+2q^2+2q^3+2q^4+4q^5+6q^6+8q^7+10q^8+14q^9+18q^{10}+\cdots .$$

\section{Correspondence with Overpartitions}
\label{Sec:3}
There is a one-to-one correspondence between 2-colored Rogers-Ramanujan type partitions and previously defined overpartitions $D_{k,a}(n)$ for k=a=2.

Let $\pi=(y_1,\cdots,y_i,y_{i+1},\cdots,y_m)$ be an arbitrary partition of $n$, first of all, for both cases all parts are distinct. Secondly, for the case that there are $t$ number of consecutive parts, for the colored case, there are only two possibilities, they should be alternatively red and black, e.g. for three consecutive parts $i$, $i+1$ and $i+2$ cases are
$$\cdots,i,\color{red}{i+1} \color{black},i+2, \cdots \;\; \text{and} \;\; \cdots, \color{red}i \color{black}, {i+1},\color{red}{i+2} \color{black},\cdots$$ 
which means two consecutive parts can not be of the same color, and for overpartition case, the first $t-1$ parts should be overlined and there are two possibilities for the last one, e.g. for three consecutive parts we will have
$$\cdots,\overline{i},\overline{i+1},i+2,\cdots \;\; \text{and} \;\; \cdots,\overline{i},\overline{i+1},\overline{i+2},\cdots$$
implies the first one should be overlined and the second one can be overlined or non-overlined. If $y_{i+1}-y_i > 1$, then there are four cases for both colored cases and overpartition one,
$$\cdots,y_i,y_{i+1}, \cdots \;\; \text{and} \;\; \cdots, y_i,\color{red} y_{i+1} \color{black},\cdots  \;\; \text{and} \;\; \cdots, \color{red}y_i \color{black}, y_{i+1},\cdots  \;\; \text{and} \;\; \cdots, \color{red}y_i \color{black}, \color{red}y_{i+1}\color{black},\cdots$$
and
$$\cdots,y_i,y_{i+1},\cdots \;\; \text{and} \;\; \cdots,y_i,\bar{y}_{i+1},\cdots  \;\; \text{and} \;\; \cdots,\bar{y}_i,y_{i+1},\cdots  \;\; \text{and} \;\; \cdots,\bar{y}_i,\bar{y}_{i+1},\cdots$$
so, there exists a one-to-one correspondence between them.

It is not hard to see another correspondence between $D_{2,1}(n)$ and the following partition type,

\begin{definition}\label{R3}
	Let $R_3(n)$ denote the number of 2-colored Rogers-Ramanujan partitions which do not allow to have a red $1$ in the partition.
\end{definition}

Then
\begin{theorem}
	For definition \ref{R3} the following identity holds
	$$R_3(n)=D_{2,1}(n)=\frac{(-q)_\infty(q^1,q^3,q^4;q^4)_\infty}{(q)_\infty}$$
\end{theorem}

\section{Future Works}
\label{Sec:4}
There are two options for further work on this topic, first one is to extend the number of colors to more than two, ideally for arbitrary number of colors. Our second suggestion is two extend it to 2-colored Gordon-Rogers-Ramanujan type partitions, i.e. for 2-colored partition $y_1+ \cdots + y_m$, and for $k \geq 1$, $y_j - y_{j+k-1} \geq 1$ if $y_j$ is black and $y_j - y_{j+k-1} \geq 2$ otherwise.

\end{document}